\documentclass[12pt]{article}
\usepackage{amsmath,bm}
\usepackage{t1enc}
\usepackage[latin1]{inputenc}
\usepackage[english]{babel}
\usepackage{amsmath,amsthm}
\usepackage{amsfonts}
\usepackage{latexsym}
\usepackage{graphicx}
\usepackage{epstopdf}
\usepackage{float}
\usepackage{tablists}
\usepackage{amssymb}
\usepackage[natural]{xcolor}
\usepackage[justification=centering]{caption}
\newtheorem{theorem}{Theorem}
\newtheorem{lemma}[theorem]{Lemma}
\newtheorem{claim}[theorem]{Claim}

\newtheorem{conjecture}[theorem]{Conjecture}

\newtheorem{corollary}[theorem]{Corollary}
\newcommand\ch{{ch}}

\title{\LARGE The total coloring of $K_5$-minor-free graphs \thanks{This work is supported by NSFC(11971270, 11631014) of China and Shandong Province Natural Science Foundation (ZR2018MA001, ZR2019MA047) of China }}
\author{Fan Yang and Jianliang Wu\thanks{Corresponding author. E-mail address: jlwu@sdu.edu.cn.}\\
{\small School of Mathematics, Shandong University, Jinan 250100, China}
 }
\date{}

\newcommand{\send}[2]{c(#1 \rightarrow #2)}

\begin{document}

\maketitle
\begin{abstract}
A total $k$-coloring of a graph $G$ is a coloring of $V(G)\cup E(G)$ using $k$ colors such that no two adjacent or incident elements receive the same color. The total chromatic number $\chi''(G)$ of $G$ is the smallest integer $k$ such that $G$ has a total $k$-coloring.  In the paper, it is proved that for any $K_5$-minor-free graph $G$, $\chi''(G)\leq \Delta(G)+2$ if $\Delta(G)\geq 7$. Moreover, $\chi''(G)=\Delta(G)+1$ if $\Delta(G)\geq 10$.

\end{abstract}

\baselineskip 0.7cm
\section{Introduction}\label{introduc}
All graphs considered in this paper are simple, finite and undirected, and we follow \cite{RD} for the terminologies and notation not defined here. Let $G=(V,E)$ be a graph. If $uv\in E(G)$, then $u$ is called a \emph{neighbor} of $v$. The \emph{degree} $d_G(v)$ (or simply $d(v)$) of a vertex $v$ is the number of neighbors of $v$. We use $\Delta(G)$ and $\delta(G)$ (or simply $\Delta$ and $\delta$) to denote the maximum degree and the minimum degree of $G$, respectively. A $k$-, $k^+$- or $k^-$-vertex is a vertex of degree $k$,  at least $k$ or at most $k$, respectively. For $v\in V(G)$, denote by $N_G(v)$ (or simply $N(v)$) the set of neighbors of $v$, and by $N_k(v)$ ($N_k^+(v)$ or  $N_k^-(v))$ the set of neighbors of $v$ of degree $k$ (at least $k$ or at most $k$, respectively). Let $|N_k(v)|=n_k(v)$, $|N_k^+(v)|=n_k^+(v)$ and $|N_k^-(v)|=n_k^-(v)$.

A \emph{total $k$-coloring} of a graph $G$ is a coloring of $V\cup{E}$ using $k$ colors such that no two adjacent or incident elements receive the same color.  The \emph{total chromatic number} $\chi''(G)$ of $G$ is the smallest integer $k$ such that $G$ has a total $k$-coloring. Behzad \cite{Be} and Vizing \cite{Viz} posed independently the following famous conjecture, which is known as the total coloring conjecture (TCC).

\vspace{2mm}{\bf Conjecture A}. For every graph $G$, $\Delta+1\leq\chi''(G)\leq\Delta+2$.

\vspace{2mm}\noindent Clearly, the lower bound is trivial. The upper bound remains open. This conjecture was confirmed for general graphs with $\Delta\leq5$. For its history, readers can see \cite{yap}. For planar graphs, the only open case is $\Delta=6$ (see \cite{Ko, SZ}). Interestingly, every planar graph with maximum degree $\Delta\geq9$ has a total $(\Delta+1)$-coloring. The result was first established in   \cite{Bo}  for $\Delta\geq14$, and was extended to $\Delta\geq12$ \cite{BWK1}, $\Delta\geq11$ \cite{BWK2}, $\Delta\geq10$ \cite{WWF}, and finally $\Delta\geq9$ \cite{KSS}.

To \emph{identify} nonadjacent vertices $x$ and $y$ of a graph $G$ is to replace these vertices by a single vertex incident with all the edges which were incident in $G$ with either $x$ or $y$. Let $e=xy$ be an edge of a graph $G$. To \emph{contract} $e$ of a graph $G$ is to delete $e$ first and then identify $x$ with $y$, finally delete one edge of any pair of parallel edges so formed. A graph $H$ is a \emph{minor} of a graph $G$ if $G$ has a subgraph contractible to $H$. $G$ is called $H$-\emph{minor-free} if $G$ does not have $H$ as a minor. In the paper we consider the total coloring of $K_5$-minor-free graphs and get the following theorem.

\begin{theorem}\label{th1}
Let $G$ be a $K_5$-minor-free graph.
\begin{itemize}
  \item[$(1)$] If $\Delta(G)\geq 7$, then $\chi''(G)\le \Delta(G)+2$;
  \item[$(2)$] If $\Delta(G)\geq 10$, then $\chi''(G)=\Delta(G)+1$.
\end{itemize}
\end{theorem}

It is well-known that a graph is planar if and only if it is $K_5$-minor-free and $K_{3,3}$-minor-free (see Theorem 10.32 in \cite{BM}). So we generalize the result of Sanders and Zhao in \cite{SZ} and that of Wang in \cite{WWF}. To prove Theorem \ref{th1}(1), we will prove a structural property of $K_5$-minor-free graphs which is similar to that of planar graphs in \cite{SZ} (see Lemma \ref{Lem-k5-deg7} in the next section). The property mentions several reducible subgraphs, all of which have been proved reducible elsewhere (see \cite{Bo,SZ}). To prove Theorem \ref{th1}(2), we also show a structural property of planar graphs (see Lemma \ref{planar}) firstly. Then we obtain some reducible subgraphs of $K_5$-minor-free graphs (see Lemma \ref{Lem-k5-deg10}). Finally, in Section \ref{section4}, we complete the proof of Theorem \ref{th1}.

\section{Structural properties of $K_5$-minor-free graphs}\label{section2}

Before proceeding, we introduce the following notation and definitions. Let $G$ be a planar graph drawn in the plane, and let $F$ be the face set of $G$. For a face $f$ of $G$, the \emph{degree $d_G(f)$} (or simply $d(f)$) is the number of edges incident with it, where each cut-edge is counted twice. Note that if there is a cut-edge, then the number of vertices incident with $f$ is less than $d_G(f)$. A \emph{$k$-face}, a \emph{$k^{-}$-face} or a \emph{$k^{+}$-face} is a face of degree $k$, at most $k$ or at least $k$, respectively. We denote a $3$-face by $f=(u, v, w)$ where $u, v, w$ are three vertices incident with $f$.  Denote by $n_d(f)$ the number of $d$-vertices incident with the face $f$, and $f_{d}(v)$ the number of $d$-faces incident with $v$. 

Note that in all figures throughout the paper, a vertex marked with $\bullet$ has no edge of $G$ incident with it other than those shown, while $\circ$ represents a vertex with uncertain degree, unless stated otherwise.

\begin{figure}[htbp]
 \begin{center}
   \includegraphics[scale=0.6]{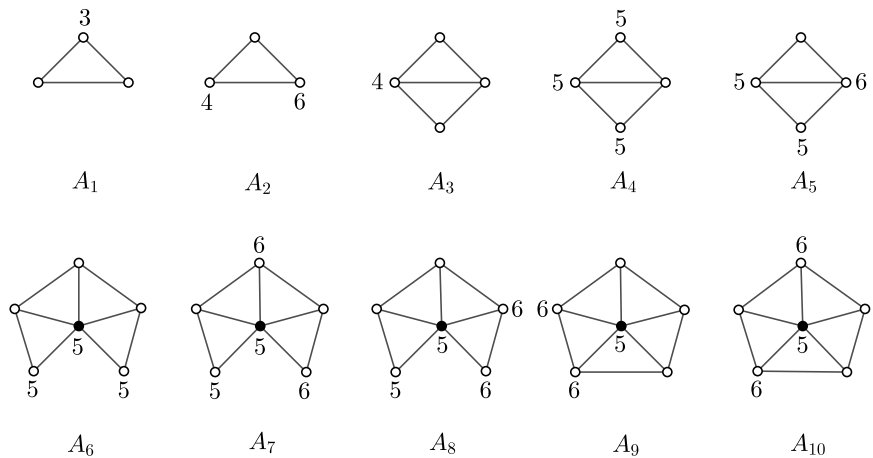}\\
   \caption[14cm]{Configurations for Lemma \ref{planar-deg7}. The number near a vertex denotes its degree and all of these numbered vertices belong to $H$.}
   \label{Fig.planar7}
 \end{center}
\end{figure}

\begin{lemma}\label{planar-deg7}
Let $G$ be a planar graph with $\delta(G)\geq1$, and let $N$ $(1\leq |N|\leq 3)$ be a set of nonadjacent vertices on the same face $f_0$ such that the subgraph $H=G-N$ has at least one edge. Suppose that
\begin{description}
    \item[$$(a)$$] $d_{G}(v)\geq3$ if $v\in V(H)$, and
    \item[$$(b)$$] $d_G(u)+ d_G(v)\geq 10$ for every edge $uv\in E(H)$.
\end{description}
Then $G$ has a subgraph isomorphic to one of the configurations $A_1$-$A_{10}$ in Figure $\ref{Fig.planar7}$.
\end{lemma}

\begin{proof}
Let $G$ be a counterexample to the lemma. Euler's formula $|V(G)|-|E(G)|+|F(G)|\geq 2$ can be expressed in the form
 $$\sum\limits_{v\in V(G)}(d_G(v)-4)+\sum\limits_{f\in F(G)}(d_G(f)-4)\leq -8. $$
It follows that
 $$\sum\limits_{v\in V(G)-N}(d_G(v)-4)+\sum_{v\in N}(d_G(v)-4+\frac{8}{|N|})+\sum_{f\in F(G)}(d_G(f)-4)\leq 0. \eqno{(\mathrm{I})}$$
An \emph{initial} charge $\ch_0$ on $V(G)\cup F(G)$ is defined by letting
\begin{align*}
\begin{split}
\ch_0(x)=\left\{
 \begin{array}{ll}
  d_G(x)-4    & $if$ \ x\in V(G)-N, \\
  d_G(x)-4+\frac{8}{|N|}    & $if$ \ x\in N, \\
  d_G(x)-4    & $if$ \ x\in F(G).
 \end{array}
 \right.
 \end{split}
\end{align*}
We will obtain a \emph{final} charge $\ch$ from $\ch_0$ by discharging rules R1-R4 below. Since these rules merely move charges around, $\mathrm{(I)}$ gives
 $$\sum\limits_{x\in V(G)\cup F(G)}\ch(x)=\sum\limits_{x\in V(G)\cup F(G)}\ch_0(x)\leq 0. \eqno{\mathrm{(II)}}$$
In the following, $\send{x}{y}$ is used to denote the amount of charges transferred from an element $x$ to an element $y$.

{\it
\begin{description}
\item[$\bf R1.$] \rm Let $v$ be a $3$-vertex of $H$ and $u\in V(G)$ be a neighbor of $v$. Then $\send{u}{v}=\frac{1}{3}$.
\item[$\bf R2.$] Let $f=(u,v,w)$ be a $3$-face of $H$ and $d_G(u) \leq d_G(v)\leq d_G(w)$. Let $f_{uv}$, $f_{uw}$ and $f_{vw}$ be the faces neighboring $f$ across edges $uv$, $uw$ and $vw$ respectively. For brevity, we write $\send{u,v,w}{f}$ to denote the triple of numbers $(\send{u}{f}, \send{v}{f}, \send{w}{f})$.
  \begin{description}
    \item[$\bf R2.1.$] If $d_G(u)\leq 4$, then $\send{u,v,w}{f}=(0, \frac{1}{2}, \frac{1}{2})$.
    \item[$\bf R2.2.$] If $d_G(u)\geq 6$ or $d_G(u)=d_G(v)=d_G(w)=5$, then $\send{u,v,w}{f}=(\frac{1}{3}, \frac{1}{3}, \frac{1}{3})$.
    \item[$\bf R2.3.$] If $d_G(u)=d_G(v)=5$ and $d_G(w)=6$, then $\send{u,v,w}{f}=(\frac{1}{6}, \frac{1}{6}, \frac{2}{3})$.
    \item[$\bf R2.4.$] If $d_G(u)=d_G(v)=5$ and $d_G(w)\ge 7$, then
\begin{align*}
\begin{split}
\send{u,v,w}{f}=\left\{
 \begin{array}{ll}
  (\frac{3}{14}, \frac{3}{14}, \frac{4}{7})    & $if$ \ d_G(f_{uw})\geq4 \ $or$ \ d_G(f_{vw})\geq4, \\
  (\frac{2}{7}, \frac{2}{7}, \frac{3}{7})    & $otherwise$.
 \end{array}
 \right.
 \end{split}
\end{align*}
   \item[$\bf R2.5.$] If $d_G(u)=5$ and $d_G(v)=d_G(w)=6$, then
   \begin{align*}
\begin{split}
\send{u,v,w}{f}=\left\{
 \begin{array}{ll}
  (\frac{1}{3}, \frac{1}{3}, \frac{1}{3})    & $if$ \ d_G(f_{uv})=d_G(f_{uw})=3, \\
  (\frac{1}{6}, \frac{1}{3}, \frac{1}{2})    & $if$ \ d_G(f_{uv})=3 \ $and$ \ d_G(f_{uw})\geq4, \\
  (\frac{1}{6}, \frac{1}{2}, \frac{1}{3})    & $if$ \ d_G(f_{uv})\geq4 \ $and$ \ d_G(f_{uw})=3, \\
  (\frac{1}{6}, \frac{5}{12}, \frac{5}{12})    & $if$ \ d_G(f_{uv})\geq4 \ $and$ \ d_G(f_{uw})\geq4. \\
 \end{array}
 \right.
 \end{split}
\end{align*}
    \item[$\bf R2.6.$] If $d_G(u)=5$, $d_G(v)=6$ and $d_G(w)\geq 7$, then
     \begin{align*}
\begin{split}
\send{u,v,w}{f}=\left\{
 \begin{array}{ll}
  (\frac{1}{14}, \frac{2}{7}, \frac{9}{14})    & $if$ \ d_G(f_{uw})\geq4, \\
  (\frac{1}{14}, \frac{1}{2}, \frac{3}{7})    & $if$ \ d_G(f_{uw})=3 \ $and$ \ d_G(f_{uv})\geq4, \\
  (\frac{5}{21}, \frac{1}{3}, \frac{3}{7})    & $otherwise$.
 \end{array}
 \right.
 \end{split}
\end{align*}
    \item[$\bf R2.7.$] If $d_G(u)=5$ and $d_G(w)\geq d_G(v)\geq 7$, then $\send{u,v,w}{f}=(\frac{1}{7}, \frac{3}{7}, \frac{3}{7})$.
  \end{description}
\item[$\bf R3.$] Let $f=(u,v,w)$ be a $3$-face of $G$ such that $w\in N$. Then  $\send{w}{f}=1$.
\item[$\bf R4.$] Let $v\in N$. If $d_G(v)=1$, then $\send{f_0}{v}=1$. If $d_G(v)\geq2$ and $d_G(f_0)\geq6$, then $\send{f_0}{v}=\frac{1}{2}$.
\end{description}}

We will now show that $\ch(x)\geq0$ for each $x\in V(G)\cup F(G)$ and $ch(v)>0$ for each $v\in N$ to get a contradiction to $\mathrm{(II)}$.

\textbf{Let $x=f$ be a face of $G$.} Then $\ch_0(f)=d_G(f)-4$. We consider two cases.

\noindent
\textbf{Case F0:} $f=f_0$.

Let $n_1$ be the number of $1$-vertices in $N$ and $n_2=|N|-n_1$. If $d_G(f_0)\geq7$, then $\ch(f_0)\geq d_G(f_0)-4-3\times 1\geq0$ by R4. If $d_G(f_0)=6$, then $n_1\leq1$ and it follows that $\ch(f_0)\geq6-4-1-2\times \frac{1}{2}=0$. If $d_G(f_0)=5$, then $n_1\leq1$ and $\ch(f_0)=5-4-1\geq0$. If $d_G(f_0)=4$, then $n_1=0$ and $\ch(f_0)=4-4=0$. If $d_G(f_0)=3$, then $n_1=0$ and $\ch(f_0)=\ch_0(f_0)+1=0$ by R3.

\noindent
\textbf{Case F1:} $f\neq f_0$.

If $d_G(f)\geq4$ then $\ch(f)=\ch_0(f)=d_G(f)-4\geq0$. Otherwise $d_G(f)=3$ and $\ch(f)\geq \ch_0(f)+1=0$ by R2 and R3.

\textbf{Now let $v=x$ be a vertex of $G$.} There are several cases.

\noindent
\textbf{Case V2:} $v\in N$.

Then $\ch(v)\geq \ch_0(v)-f_3(v)-\frac{n_3(v)}{3}=\frac{8}{|N|}-4+d_G(v)-f_3(v)-\frac{n_3(v)}{3}$ by R1 and R3. Since $G$ does not contain $A_1$, any $3$-neighbor of $v$ is not incident with a $3$-face of $G$. If $d_G(f_0)=3$, then $|N|=1$ and $f_3(v)+\frac{n_3(v)}{3}\leq f_3(v)+n_3(v)\leq d_G(v)$, and it follows that $\ch(v)\geq \frac{8}{|N|}-4>0$. So assume $d_G(f_0)\geq4$. If $d_G(v)=1$, then $f_3(v)+\frac{n_3(v)}{3}\leq \frac{1}{3}< d_G(v)$, otherwise $f_3(v)+\frac{n_3(v)}{3}\leq d_G(v)-1$. If $|N|\leq2$, then $\ch(v)>\frac{8}{|N|}-4\geq0$, and thus $\ch(v)>0$. So assume $|N|=3$. Since any two vertices in $N$ are not adjacent, we get $d_G(f_0)\geq6$. If $d_G(v)\geq2$, then $\send{f_0}{v}=\frac{1}{2}$ by R4 and $\ch(v)\geq \frac{8}{3}-4+1+\frac{1}{2}>0$. If $d_G(v)=1$, then $\send{f_0}{v}=1$ by R4 and $\ch(v)\geq \frac{8}{3}-4+1+1-\frac{1}{3}>0$.

\noindent
\textbf{Case V3:} $v\in V(H)$ and $d_G(v)=3$.

Then $\ch_0(v)=-1$ and $\ch(v)=-1+3\times \frac{1}{3}=0$ by R1.

\noindent
\textbf{Case V4:} $v\in V(H)$ and $d_G(v)=4$.

Then $\ch(v)=\ch_0(v)=0$.

\noindent
\textbf{Case V5:} $v\in V(H)$ and $d_G(v)=5$.

Then $\ch_0(v)=1$. By the hypothesis (b), every neighbor of $v$ in $H$ is a $5^{+}$-vertex, and so $v$ sends no charges to its neighbors by R1. According to R2 and R3, $v$ sends at most $\frac{1}{3}$ to each incident $3$-face, but nothing if the $3$-face has a vertex in $N$. If $v$ has a neighbor in $N$, then $v$ sends charges to at most three $3$-faces and $\ch(v)\geq1-\frac{3}{3}=0$. So we may assume that all neighbors of $v$ are in $H$. If $f_3(v)\leq3$, then again $\ch(v)\geq1-\frac{3}{3}=0$, so we may assume that $f_3(v)\geq4$. We consider two subcases separately.

\noindent
\textbf{Subcase V5.1:} $f_3(v)=4$.

Let the four 3-faces incident with $v$ be $f_1=(v, v_1, v_2)$, $f_2=(v, v_2, v_3)$, $f_3=(v, v_3,v_4)$ and $f_4=(v, v_4, v_5)$. W.l.o.g., we assume $d_G(v_1)\leq d_G(v_5)$. By the absence of $A_6$, the possible values for $(d_G(v_1), d_G(v_5))$ are $(5,6)$, $(5,7^{+})$ and $(6^{+},6^{+})$. Furthermore, by the absence of $A_4$-$A_8$, the possibilities of degree sequences $v_1, \ldots, v_5 $ as well as the maximum total charges $\Sigma_4$ from $v$ to four faces $f_1, \ldots, f_4$ are tabulated as follows.
\begin{center}
$
\begin{array}{llllllllll}
d_G(v_1)&d_G(v_2)&d_G(v_3)&d_G(v_4)&d_G(v_5) \\
5&7^{+}&7^{+}&7^{+}&6 \\
5&7^{+}&6&6^{+}&7^{+} \\
5&7^{+}&7^{+}&5^{+}&7^{+} \\
6^{+}&5&7^{+}&5^{+}&6^{+} \\
6^{+}&6&6&6^{+}&6^{+} \\
6^{+}&6&7^{+}&5^{+}&6^{+} \\
6^{+}&7^{+}&5&7^{+}&6^{+} \\
6^{+}&7^{+}&6&6^{+}&6^{+} \\
6^{+}&7^{+}&7^{+}&5^{+}&6^{+}
\end{array}
$
\ \ \ \ \ \
$
\begin{array}{llllllllll}
\Sigma_4 \\
\frac{2}{7}+\frac{1}{7}+\frac{1}{7}+\frac{1}{14}<1, \\
\frac{2}{7}+\frac{5}{21}+\frac{1}{3}+\frac{1}{7}=1, \\
\frac{2}{7}+\frac{1}{7}+\frac{2}{7}+\frac{3}{14}<1, \\
\frac{3}{14}+\frac{2}{7}+\frac{2}{7}+\frac{3}{14}=1, \\
\frac{1}{6}+\frac{1}{3}+\frac{1}{3}+\frac{1}{6}=1, \\
\frac{1}{6}+\frac{5}{21}+\frac{2}{7}+\frac{3}{14}<1, \\
\frac{1}{7}+\frac{2}{7}+\frac{2}{7}+\frac{1}{7}<1, \\
\frac{1}{7}+\frac{5}{21}+\frac{1}{3}+\frac{1}{6}<1, \\
\frac{1}{7}+\frac{1}{7}+\frac{2}{7}+\frac{3}{14}<1.
\end{array}
$
\end{center}

In each case $\ch(v)\geq\ch_0(v)-\Sigma_4\geq0$.

\noindent
\textbf{Subcase V5.2:} $f_3(v)=5$.

Suppose first that $v$ is adjacent to a $5$-vertex $u$. Then every other neighbor of $v$ is $7^{+}$-vertex by the absence of $A_4$, $A_5$ and $A_{10}$, and it follows that $v$ sends at most $\frac{2}{7}$ to each face incident with $uv$ by R2.4, and at most $\frac{1}{7}$ to each other incident face by R2.7. So $\ch(v)\geq1-2\times\frac{2}{7}-3\times\frac{1}{7}=0$.

Suppose that all neighbors of $v$ are $6^{+}$-vertices. Then at most one of these neighbors is a $6$-vertex by the absence of $A_9$ and $A_{10}$. So $\ch(v)\geq 1- 2\times \frac{5}{21}- 3\times \frac{1}{7}>0$ by R2.6 and R2.7.

\noindent
\textbf{Case V6:} $v\in V(H)$ and $d_G(v)=6$.

Then $\ch_0(v)=2$. By the hypothesis (b), every neighbor of $v$ in $H$ is a $4^{+}$-vertex, and then $v$ sends no charges to its neighbors by R1. So $v$ sends at most $\frac{2}{3}$ to each incident $3$-face by R2. If $f_3(v)\leq3$, then $\ch(v)\geq2-3\times\frac{2}{3}=0$. If $f_3(v)=4$, then there are three possibilities shown in Figure \ref{Fig.p7-6v}. For Figure \ref{Fig.p7-6v}$(B_1)$, $v$ sends at most $\frac{1}{2}$, $\frac{1}{3}$, $\frac{1}{3}$, $\frac{1}{2}$ to $f_1$, $f_2$, $f_3$, $f_4$ respectively. For Figure \ref{Fig.p7-6v}$(B_2)$, $v$ sends at most $\frac{1}{2}$, $\frac{1}{3}$, $\frac{1}{2}$, $\frac{2}{3}$ to $f_1$, $f_2$, $f_3$, $f_5$ respectively. For Figure \ref{Fig.p7-6v}$(B_3)$, $v$ sends at most $\frac{1}{2}$, $\frac{1}{2}$, $\frac{1}{2}$, $\frac{1}{2}$ to $f_1$, $f_2$, $f_4$, $f_5$ respectively. In each case the total charges are at most $2$ and $\ch(v)\geq2-2=0$.

\begin{figure}[htbp]
 \begin{center}
   \includegraphics[scale=0.8]{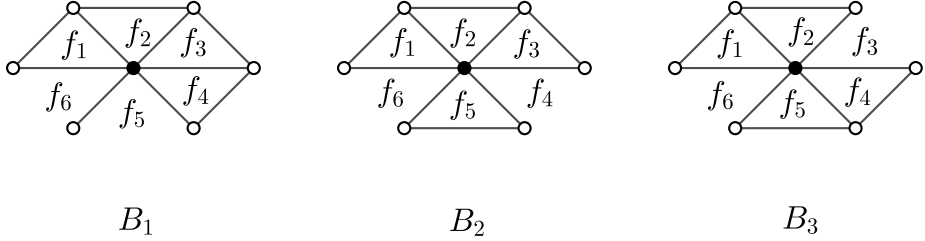}\\
   \caption{Distribution diagram of four $3$-faces of the $6$-vertex $v$.}
   \label{Fig.p7-6v}
 \end{center}
\end{figure}

If $f_3(v)=5$, then the five consecutive $3$-faces receive at most $\frac{1}{2}$, $\frac{1}{3}$, $\frac{1}{3}$, $\frac{1}{3}$, $\frac{1}{2}$, respectively. Thus $ch(v)\geq2-2\times\frac{1}{2}-3\times\frac{1}{3}=0$. If $f_3(v)=6$, then each of the six $3$-faces receives at most $\frac{1}{3}$ from $v$. Thus $\ch(v)\geq2-6\times\frac{1}{3}=0$.

\noindent
\textbf{Case V7:} $v\in V(H)$ and $d_G(v)\geq7$.

Note that for each incident 3-face $f$ of $v$, $\send{v}{f}\leq \max\{\frac13, \frac{4}{7}, \frac 37, \frac{9}{14}\} = \frac{9}{14}$ by R2. Suppose that $f_1=(v, v_1, v_2)$, $f_2=(v, v_2, v_3), \cdots,  f_r=(v, v_{r}, v_{r+1})$ be consecutive 3-faces in a clockwise order, where $r\geq 2$. Then $\max\{\send{v}{f_1}, \send{v}{f_{r}}\}\leq \frac{9}{14}$, $\send{v}{f_i}\leq \frac 37$ $(2\leq i<r)$. So
$$\send{v}{f_1}+ \cdots+ \send{v}{f_r}\leq \frac 37(r+1).\eqno{\mathrm{(III)}}$$

Suppose that $d_G(v)=7$. If $n_3(v)\geq 6$, then $f_3(v)=0$ by $A_1$ and $ch(v)\geq3- \frac13 \times7= \frac23>0$. If $3\leq n_3(v)\leq 5$, then $f_3(v)+ n_3(v) \leq 6$ and it follows from R2 that $ch(v)\geq 3- \frac13 n_3(v) - \frac{9}{14} f_3(v) \geq \frac{1}{14}>0$. Suppose $n_3(v)=2$. Then $f_3(v)\leq 4$. If $f_3(v)\leq 3$, then $ch(v) \geq 3- \frac13\times 2- \frac{9}{14} \times 3>0$. Otherwise $f_3(v)=4$ and the four 3-faces must be consecutive, it follows from $\mathrm{(III)}$ that $ch(v)\geq 3-\frac13 \times 2-\frac{3}{7}\times5>0$. Suppose $n_3(v)=1$. Then $f_3(v)\leq5$. If $f_3(v)\leq 4$, then $ch(v) \geq 3- \frac13- \frac{9}{14} \times 4>0$. Otherwise $f_3(v)=5$ and it follows from $\mathrm{(III)}$ that $ch(v)\geq 3-\frac13-\frac{3}{7}\times 6>0$. Hence we assume that $n_3(v)=0$.

If $f_3(v)\leq 4$, then $ch(v)\geq 3-\frac{9}{14}\times 4>0$. If $f_3(v)=5$, then $ch(v)\geq3-\max\{\frac{9}{14}\times2+\frac{3}{7}\times3, \frac{9}{14}\times3+\frac{3}{7}\times2, \frac{9}{14}\times4+\frac{3}{7}\}=0$. If $f_3(v)=6$, then it follows from $\mathrm{(III)}$ that $ch(v) \geq 3- \frac{3}{7}\times 7=0$. Suppose that $f_3(v)=7$. Then $v$ sends at most $\frac37$ to each of its incident 3-face by R2. Thus $ch(v)\geq 3- \frac37 \times 7=0$.

Suppose that $d_G(v)=8$ and $n_3(v)\geq 1$. If $n_3(v)\geq 7$, then $f_3(v)=0$ by $A_1$ and $ch(v)\geq4- \frac13 \times8= \frac43>0$. If $2\leq n_3(v)\leq 6$, then $f_3(v)+ n_3(v) \leq 7$ and it follows from R2 that $ch(v)\geq 4- \frac13 n_3(v) - \frac{9}{14} f_3(v)>0$. Suppose $n_3(v)=1$. Then $f_3(v)\leq6$. If $f_3(v)\leq 5$, then $ch(v) \geq 4- \frac13- \frac{9}{14} \times 5>0$. Otherwise $f_3(v)=6$ and it follows from $\mathrm{(III)}$ that $ch(v)\geq 4-\frac13-\frac{3}{7}\times 7>0$. Hence we assume that $n_3(v)=0$, this case will be considered later.

Suppose that $d_G(v)\geq 9$ and $n_3(v)\geq 1$. If $n_3(v)\geq d_G(v)-1$, then $f_3(v)=0$ by $A_1$ and $ch(v)\geq d_G(v)-4 - \frac13 d_G(v)>0$. If $1\leq n_3(v)\leq d_G(v)-2$, then $f_3(v)+ n_3(v) \leq d_G(v)-1$ and it follows from R2 that $ch(v)\geq d_G(v)- 4- \frac13 n_3(v)-\frac{9}{14} f_3(v) \geq \frac{15(d_G(v)-9)+7}{42}>0$.

Finally, suppose that $d_G(v)\geq 8$ and $n_3(v)=0$. If $f_3(v)\leq d_G(v)-3$, then $ch(v)\geq d_G(v)- 4- \frac{9}{14} (d_G(v)-3)=\frac{5(d_G(v)-7)+6}{14}>0$. Suppose that $f_3(v)= d_G(v)-2$. If the $3$-faces incident with $v$ are not all consecutive, but are in two consecutive runs of lengths $r_1$ and $r_2$ ($r_1+r_2=d_G(v)-2$), then it follows from $\mathrm{(III)}$ that $ch(v) \geq d_G(v)-4-\frac37 (r_1+1)-\frac37 (r_2+1)=\frac{4}{7}d_G(v)-4>0$. Otherwise, the $3$-faces incident with $v$ are all consecutive and it follows from $\mathrm{(III)}$ that $ch(v) \geq d_G(v)- 4- \frac{3}{7}(d_G(v)-2+1)>0$. If $f_3(v)\leq d_G(v)-1$, we also have $ch(v) \geq d_G(v)- 4- \frac{3}{7}(d_G(v)-1+1) =\frac{4}{7}d_G(v)-4>0$. If $f_3(v)= d_G(v)$, then $ch(v) \geq d_G(v)- 4- \frac37 d_G(v) = \frac{4}{7}d_G(v)-4>0$.

This completes the proof of Lemma \ref{planar-deg7}.
\end{proof}

Let $G$ be a connected graph, $T$ be a tree, and $ \mathcal{F}=\{V_t\subset V(G):t\in V(T)\}$ be a family such that $V_i\setminus V_j\neq \emptyset$ for any two different vertices $i$ and $j$ in $T$. The ordered set $(T, \mathcal{F})$ is called a \emph{tree-decomposition} of $G$ if it satisfies the following three conditions:
\begin{description}
\item[(T1)] $V(G)=\bigcup_{t\in V(T)}V_t$,

\item[(T2)] for every edge $e\in E(G)$, there exists a $t\in V(T)$ such that both ends of $e$ lie in $V_t$,

\item[(T3)] if $t_1$, $t_2$, $t_3\in V(T)$ and $t_2$ is on the $(t_1, t_3)$-path connecting $t_1$ and $t_3$ of $T$, then $V_{t_1}\cap V_{t_3}\subset V_{t_2}$.
\end{description}

Conditions (T1) and (T2) together say that $G$ is the union of the subgraphs $G[V_t]$ $(t\in V(T))$, each of these subgraphs is called a \emph{part} of $(T, \mathcal{F})$. Condition (T3) implies that the parts of $(T, \mathcal{F})$ are organized roughly like a tree. More explanations and results refer to \cite{RD}. Consider a fixed tree-decomposition $(T, \mathcal{F})$ of $G$, if there are two vertices $s$, $t\in V(T)$ such that $V_s\cap V_t\neq \emptyset$, then $V_s\cap V_t$ forms a vertex cut of $G$, called a \emph{separator} \emph{set} of $(T, \mathcal{F})$. A separator set is \emph{simplicial} if it induces a complete subgraph in $G$. If every separator set of $(T, \mathcal{F})$ is simplicial and contains at most $k$ vertices, then $(T, \mathcal{F})$ is called $k$-\emph{simplicial}.

\begin{figure}[htbp]
\begin{center}
\includegraphics[scale=0.75]{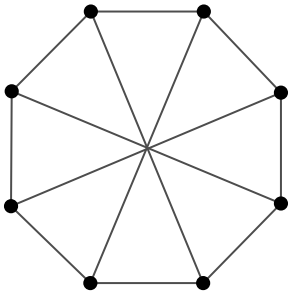} \\
\caption{The Wagner graph $W$}
\label{Wagner}
\end{center}
\end{figure}

Wagner \cite{Wagner} proved the following result about the tree-decomposition of $K_5$-minor-free graphs.
\begin{lemma}{\rm \cite{Wagner}}
Let $G$ be an edge-maximal $K_5$-minor-free graph and $|G|\geq 4$. Then $G$ has a $3$-simplicial tree-decomposition $(T, \mathcal{F})$ such that each part is a planar triangulation or the Wagner graph $W$ $($see Figure $\ref{Wagner}$$)$.
\end{lemma}

The lemma can be also found in \cite{RD} (Theorem 7.3.4). It implies that any $K_5$-minor-free graph $G$ has a supergraph $H$ of the same order with a simplicial tree-decomposition into plane triangulations and the Wagner graph $W$. Thus we regard a tree-decomposition of $H$ as that of $G$. Given a $K_5$-minor-free graph $G$, a tree-decomposition $(T, \mathcal{F})$ of $G$ is called \emph{$k$-simplified} if the resulting graph obtained from $G$ is still $K_5$-minor-free by adding extra edges as necessary such that every separator set is simplicial.  So the following corollary holds.

\begin{corollary}\label{wagner}
Let $G$ be a $K_5$-minor-free graph and $|G|\geq 4$. Then $G$ has a $3$-simplified tree-decomposition $(T, \mathcal{F})$ such that each part is a planar graph or the Wagner graph $W$.
\end{corollary}

In other words, $G$ can be composed of $G_1, \ldots, G_p$ by repeatedly taking $1$-, $2$-, and $3$-sums, where a \emph{$k$-sum} of two graphs $G_1$ and $G_2$ is obtained from their disjoint union $G_1 \cup G_2$ by identifying a $k$-clique of $G_1$ with a $k$-clique of $G_2$, and then possibly deleting some edges of this common clique (see \cite{Ding} or Exercises 10.5.10 of \cite{BM}).

For a fixed tree-decomposition $(T,\mathcal{F})$, we always choose a leaf of $T$ as the root to obtain a rooted tree, also denoted by $T$. Let $ab$ be an edge of $T$. If $b$ is closer to the root than $a$, then $a$ is called a \emph{child} of $b$ and $b$ is the \emph{parent} of $a$. We denote $S_a=V_a\cap V_b$, $G'_a=G[V_a]-S_a$, and $G^*_a=G[V_a]-E(G[S_a])$.

\begin{lemma}\label{Lem-k5-deg7}
Let $G$ be a connected $K_5$-minor-free graph. Then one of the following conditions holds.
\begin{description}
  \item[$(1)$] $\delta(G)\leq 2$,
  \item[$(2)$] $G$ contains an edge $uv$ such that $d(u)+d(v)\leq 9$,
  \item[$(3)$] $G$ contains two $3$-vertices $u$ and $v$ satisfying $|N(u)\cap N(v)|\geq 2$,
  \item[$(4)$] $G$ contains a subgraph isomorphic to one of the configurations in Figure $\ref{Fig.planar7}$.
\end{description}
\end{lemma}
\begin{proof} Let $G$ be a counterexample to the lemma  with $|V|+|E|$ as small as possible. Then
\begin{description}
  \item[$(\mathrm{a})$] $\delta(G)\geq 3$, and
  \item[$(\mathrm{b})$] $d(u)+d(v)\geq 10$ for every edge $uv\in E(G)$.
\end{description}

By Corollary \ref{wagner}, $G$ has a $3$-simplified tree-decomposition $(T, \mathcal{F})$ such that each part is a planar graph or the Wagner graph, and assume that $T$ contains as many vertices as possible and subject to that the number of leaves is as large as possible. Suppose that $|V(T)|=1$. If $G$ is the Wagner graph, then $d_G(u)+d_G(v)=6$ for every edge $uv\in E(G)$, a contradiction to $(\mathrm{b})$. If $G$ is a planar graph, then it follows from Lemma \ref{planar-deg7} that $G$ contains one of the configurations in Figure $\ref{Fig.planar7}$ (we can arbitrarily choose a vertex of $G$ as an element of $N$ in Lemma \ref{planar-deg7}), which is also a contradiction. So $|V(T)|\geq 2$ and it implies that $T$ has at least one leaf. The following claim is heavily used in completing the proof.

\begin{claim}\label{cl0}
For every leaf $x$ of $T$, we have that $|S_x|=3$ and  $G^*_x=K_{1,3}$ is a star of order $4$.
\end{claim}

\begin{proof}
Let $a$ be a leaf of $T$. By $(\mathrm{b})$, $G^*_a$ is not the Wagner graph. So it is planar, and it is connected by the maximality of $|V(T)|$. If $|V(G'_a)|\geq 2$, then $(4)$ holds by Lemma \ref{planar-deg7} (here $G'_a$, $S_a$ plays the role of $H$, $N$ in Lemma \ref{planar-deg7}, respectively), a contradiction. If $|V(G'_a)|=1$ and $|S_a|\leq 2$, then $\delta(G)\leq 2$, a contradiction to $(\mathrm{a})$. Hence $|S_a|=3$ and $G^*_a=K_{1,3}$ is a star of order $4$.
\end{proof}

Suppose that $|V(T)|=2$, and let $V(T)=\{v_1,v_2\}$. Suppose that  $v_1$ is the root of $T$. Then $v_2$ is the leaf of $T$. It follows from Claim \ref{cl0} that $G^*_{v_2}=K_{1,3}$. In fact, we can also choose $v_2$ as the root of $T$. Thus $G^*_{v_1}=K_{1,3}$, too. Thus $|V(G)|=5$, a contradiction to $(\mathrm{b})$. So $|V(T)|\geq 3$.

Now we choose a vertex $c$ of $T$, $N_T(c)=\{c_0, c_1, ..., c_t\}$ $(t\geq 1)$, such that $c_0$ is the parent of $c$ and $c_1, ..., c_t$ are leaves of $T$. By Claim \ref{cl0}, we have that for each $i\in \{1,\ldots,t\}$,  $|S_{c_i}|=3$ and  $G^*_{c_i}$ is a star of order $4$. Thus, we denote $V(G^*_{c_i})=\{u_i, u_{i1}, u_{i2}, u_{i3}\}$ and $S_{c_i}=\{u_{i1}, u_{i2}, u_{i3}\}$ $(1\leq i\leq t)$.  Let $G^{+}$ be the graph obtained from $G$ by adding extra edges as necessary such that every separator set is simplicial, that is, induces a complete subgraph of $G^{+}$. By our choice of tree-decomposition, $G^{+}$ is also a $K_5$-minor-free graph, and every part of $(T, \mathcal{F})$ in $G^+$ is planar or the Wagner graph. Note that $G^{+}[V_c]$ is not the Wagner graph since $S_{c_1}=\{u_{11}, u_{12}, u_{13}\}\subseteq V_c$ and it induces a triangle in $G^{+}$, whereas the Wagner graph is triangle-free. So $G^{+}[V_c]$ is planar and then we embed it in the plane so that the vertices and edges of the complete graph $G^{+}[S_c]$ are all in the boundary of the outside face. By the maximality of $|V(T)|$, $G^{+}[V_c]$ has no separating triangle, as what is inside the triangle would imply the existence of a neighbor of $c$ in $T$ that is different from $c_0,c_1,\ldots,c_t$. Thus, for each $j$ ($1\leq j\leq t$),  $u_{j1}, u_{j2}, u_{j3}$ are the three incident vertices of some $3$-face $f_j$ of $G^{+}[V_c]$. By (3), $|S_{c_i}\cap S_{c_j}|\leq1$ for any $i,j$ $(1\leq i<j\leq t$). So we can embed every star $G^{*}_{c_j}$ into $f_j$ ($1\leq j\leq t$) to form a new planar graph $K^{+}$.

Now let $K=G^{*}_c\cup G_{c_1}\cup \cdots \cup G_{c_t}$ and $\overline{K}=K-S_c$. Since $K\subseteq K^{+}$, $K$ is planar. By (a), every vertex of $\overline{K}$ has degree at least $3$ in $K$. Suppose that for every $i$ $(1\leq i\leq t)$, $S_{c_i}\subseteq S_c$. Then we could change $T$ by making $c_1, \ldots, c_t$ adjacent to $c_0$ instead of to $c$, and this would give a valid tree-decomposition with a larger number of leaves, which contradicts to our choice of $T$. So there exists some $i$ $(1\leq i\leq t)$ such that $S_{c_i}\not\subseteq S_c$, it implies that there is at least one edge $uv\in E(G_{c_i})$ such that $u$, $v\notin S_c$. Hence, $\overline{K}$ has at least one edge. By Lemma \ref{planar-deg7} (where $K, \overline{K}, S_c$ plays the role of $G, H, N$ in Lemma \ref{planar-deg7}, respectively), (4) holds for $K$ and then for $G$. This contradiction completes the proof of the lemma.
\end{proof}

In the following, we prove another structural property of planar graphs.

\begin{figure}[htbp]
 \begin{center}
   \includegraphics[scale=0.8]{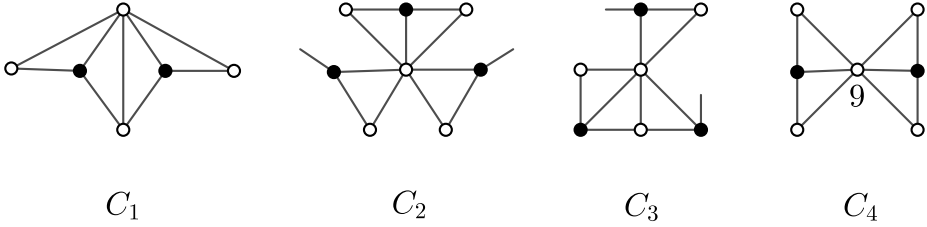}\\
   \caption{Configurations for Lemma \ref{planar}. The number near a vertex in $C_4$ denotes its degree, and all vertices marked with $\bullet$ as well as the marked $9$-vertex in $C_4$ belong to $H$.}
   \label{Fig.planar10}
 \end{center}
\end{figure}

\begin{lemma}\label{planar}
Let $G$ be a planar graph with $\delta(G)\geq1$, and let $N$ $(1\leq |N|\leq 3)$ be a set of nonadjacent vertices on the same face $f_0$ such that the subgraph $H=G-N$ has at least one edge. Suppose that
\begin{description}
    \item[$$(a)$$] $d_G(v)\geq 3$ for each vertex $v\in V(H)$, and
    \item[$$(b)$$] $d_G(u)+d_G(v)\geq 12$ for every edge $uv\in E(H)$.
\end{description}
Then $G$ has a subgraph isomorphic to one of the configurations $C_1$-$C_4$ in Figure $\ref{Fig.planar10}$.
\end{lemma}
\begin{proof}
Let $G$ be a counterexample to the lemma  with $|V|+|E|$ as small as possible. Then $G$ is connected, Euler's formula $|V(G)|-|E(G)|+|F(G)|\geq 2$ can be expressed in the form
 $$\sum\limits_{v\in V(G)}(d_G(v)-6)+\sum\limits_{f\in F(G)- f_0}(2d_G(f)-6)+(2d_G(f_0)+6)\leq 0. \eqno{\mathrm{(IV)}}$$
Thus an \emph{initial} charge $\ch_0$ on $V(G)\cup F(G)$ is defined as follows.
$$\ch_0(x)=\left\{
  \begin{array}{ll}
    d_G(x)-6, & \hbox{if \ $x\in V(G)$,} \\
    2d_G(x)-6, & \hbox{if \ $x\in F(G)-\{f_0\}$,} \\
    2d_G(x)+6, & \hbox{if \ $x=f_0$.}
  \end{array}
\right.$$
We will obtain a \emph{final} charge $\ch$ from $\ch_0$ by discharging rules R1-R4 below. Since these rules merely move charges around, $\mathrm{(IV)}$ gives
$$\sum_{x\in{{V(G)}\cup{F(G)}}}\ch(x)=\sum_{x\in{{V(G)}\cup{F(G)}}}\ch_0(x)\leq 0.\eqno{\mathrm{(V)}}$$
But a contradiction will be obtained by proving that $\ch(x)\geq0$ for each element $x$, with strict inequality in at least one case.

A $4^{+}$-face $f\neq f_0$ is called a \emph{lo-face} if $d_G(f)=4$ and $f$ is incident with two $5^{-}$-vertices that are in $H$, otherwise $f$ is a \emph{hi-face}. For a vertex $v$ that is incident with a face $f$, $t(f,v)$ will denote the number of times that $v$ occurs in the boundary walk of $f$ (or in the component of the boundary that contains $v$, if the boundary of $f$ is disconnected), this is the number of blocks of $G$ that contain $v$ and have at least one edge incident with $f$, which is $1$ unless $v$ is a cut vertex. We use  $\send{x}{y}$ to denote the amount of charges transferred from an element $x$ to an element $y$.

\begin{description}
  \item[\textbf{R$1$.}] Let $v$ be a vertex of $G$ incident with $f_0$.
    \begin{description}
       \item[\textbf{R$1.1$.}] Suppose $v\in N$. If $d_G(v)=1$, then $\send{f_0}{v}=5$.  Otherwise $\send{f_0}{v}=4$.
       \item[\textbf{R$1.2$.}] Suppose $v\notin N$. Then $\send{f_0}{v}=1+t(f_0,v)$.
    \end{description}
  \item[\textbf{R$2$.}] Let $v\in V(H)$ be a $5^{-}$-vertex of $G$ that is incident with a $4^{+}$-face $f\neq f_0$. If $f$ is a lo-face, then $\send{f}{v}=1$. Otherwise $\send{f}{v}=2t(f,v)$.
  \item[\textbf{R$3$.}] Let $v\in N$ and let $vv'$ be an edge that is not incident with $f_0$. Then $\send{v}{v'}=1$.
  \item[\textbf{R$4$.}] Let $v\in V(H)$ be a $5^{-}$-vertex of $G$, and let $vv'$ be an edge such that $v'\in V(H)$ and $vv'$ is not incident with $f_0$.
     \begin{description}
     \item[\textbf{R$4.1$.}] If $d_G(v)=3$, then
     \begin{align*}
\begin{split}
\send{v'}{v}=\left\{
 \begin{array}{ll}
  1    & $if$ \ vv' \ \rm{is \ incident \ with \ two} \ \rm {\operatorname{3-faces}}, \\
  \frac{1}{2}    & $if$ \ vv' \ \rm{is \ incident \ with \ a} \ \rm{\operatorname{3-face}} \ and \ a \ \rm{\operatorname{lo-face}}.
 \end{array}
 \right.
 \end{split}
\end{align*}
     \item[\textbf{R$4.2$.}] If $d_G(v)=4$, then $\send{v'}{v}=\frac{1}{2}$ if $vv'$ is incident with two $3$-faces.
     \item[\textbf{R$4.3$.}] If $d_G(v)=5$, then $\send{v'}{v}=\frac{1}{5}$ if $vv'$ is incident with two $3$-faces.
     \end{description}
\end{description}

We begin to show that $\ch(x)\geq0$ for each $x\in{{V(G)}\cup{F(G)}}$ and $ch(f_0)+\sum_{v\in V(H)}ch(v)>0$ to get a contradiction to $\mathrm{(V)}$.

\textbf{Let $x=f$ be a face of $G$.} There are two cases to be considered.

\noindent
\textbf{Case F0:} $f=f_0$.

Recall that $\ch_0(f_0)=2d_G(f_0)+6$. Let $n_1^{*}$ be the number of $1$-vertices in $N$, and let $\Sigma_0=\sum\{t(f_0,v)-1: v\in V(H), v \ $is$ \ $incident$ \ $with$ \ f_0\}$. Then $n_1^{*}\leq \Sigma_0$, and the number of distinct vertices of $f_0$ is at most $d_G(f_0)-\Sigma_0$. In the worst case, by R1, $f_0$ gives $2$ to each of its incident vertices, plus an additional $\Sigma_0$ to all incident vertices in $V(H)$, plus $2$ to each vertex in $N$, plus a further $1$ to each $1$-vertex in $N$. This is a total of at most $2(d_G(f_0)-\Sigma_0)+\Sigma_0+2|N|+n_1^{*}\leq 2d_G(f_0)+2|N|$, since $n_1^{*}\leq \Sigma_0$. So $\ch(f_0)\geq (2d_G(f_0)+6)-(2d_G(f_0)+2|N|)=6-2|N|\geq0$, since $|N|\leq3$.


\noindent
\textbf{Case F1:} $f\neq f_0$.

If $d_G(f)=3$ then $\ch(f)=\ch_0(f)=2d_G(f)-6=0$. If $d_G(f)=4$, $\ch(f)\geq (2d_G(f)-6)-\max\{2,1+1\}=0$ by R2. If $d_G(f)\geq5$, then $\ch(f)\geq (2d_G(f)-6)-2\lfloor \frac{d_G(f)}{2}\rfloor\geq0$.

\textbf{Now let $x=v$ be a vertex of $G$.} Then $\ch_0(v)=d_G(v)-6$. There are several cases.

\noindent
\textbf{Case V1: $v\in N$.}

If $d_G(v)=1$, then $\ch(v)=d_G(v)-6+5=0$ by R1. If $d_G(v)=2$, then $\ch(v)=d_G(v)-6+4=0$ by R1. If $d_G(v)\geq3$, then $v$ sends $1$ along every edge $vv'$ not incident with $f_0$ by R3, and then $\ch(v)=(d_G(v)-6)+4-(d_G(v)-2)=0$.

\noindent
\textbf{Case V2:} $v\in V(H)$ and $d_G(v)\leq6$.

By (b), $v$ is not adjacent to any $5^{-}$-vertex of $G$ in $V(H)$, and so $v$ gives out nothing by R4. If $d_G(v)=6$, then $\ch(v)= \ch_0(v)=0$. So we may assume $3\leq d_G(v) \leq5$. There are two subcases.

\noindent
\textbf{Case V2.1:} $v$ is incident with $f_0$.

Then $c(f_0\rightarrow v)=1+t(f_0,v)\geq 2$. If $t(f_0,v)\geq 2$ or $d_G(v)\geq 4$, then $ch(v)\geq d_G(v)-6+1+t(f_0,v)\geq 0$. So we assume that $d_G(v)=3$ and $t(f_0,v)=1$. Thus $v$ is incident with two faces different from $f_0$ and separated by an edge $vv'$. If one of the two faces incident with $vv'$ is a $4^{+}$-face, then $v$ gets at least an additional $1$ from it by R2, otherwise both are $3$-faces and $v$ gets $1$ from $v'$ by R4.1. In both cases we have have $ch(v)\geq (d_G(v)-6)+2+1=0$.

\noindent
\textbf{Case V2.2:} $v$ is not incident with $f_0$.

In order to add the contributions to $v$ by R2, R3 and R4, we assume that all charges given to $v$ passes along an incident edge. So if $uv$ and $vw$ are consecutive edges incident with $v$ in cyclic order, and the face between them is $f$, then $f$ sends $1$ to $v$ along each of these edges ($uv$ and $vw$) if $f$ is a hi-face, and $\frac{1}{2}$ along each of them if $f$ is a lo-face by R2. Thus $v$ receives at least $1$ along an edge $uv$ by R2, except in the following cases, (i) the faces neighboring $uv$ are a $3$-face and a lo-face, in which case $v$ receives $\frac{1}{2}$ along $uv$ by R2 and, if $d_G(v)=3$, at least a further $\frac{1}{2}$ along $uv$ by R3 or R4.1. (ii) the faces neighboring $uv$ are both $3$-faces, in which case $v$ receives nothing along $uv$ by R2 but at least $1$ or $\frac{1}{2}$ or $\frac{1}{5}$ along $uv$ by R3 or R4. Thus $v$ receives at least $1$ along each incident edge if $d_G(v)=3$, at least $\frac{1}{2}$ if $d_G(v)=4$, and at least $\frac{1}{5}$ if $d_G(v)=5$, a total of $3$, $2$ and $1$ in the three cases. Thus $\ch(v)\geq0$ in each case.

\noindent
\textbf{Case V3:} $v\in V(H)$ and $7\leq d_G(v)\leq8$.

If $d_G(v)=7$, then the neighbors of $v$ in $H$ are $5^{+}$-vertices. Since two $5$-vertices cannot be adjacent, there are at most three neighbors of $v$ that receive charges from $v$ by R4.3, so $\ch(v)\geq \ch_0(v)-3\times \frac 1 5=\frac 2 5>0$. If $d_G(v)=8$, then it sends at most $\frac 12$ to at most four neighbors of $v$ by R4.2 and R4.3, so $\ch(v)\geq \ch_0(v)-4\times \frac 12=0$.

\noindent
\textbf{Case V4:} $v\in V(H)$ and $d_G(v)\geq 9$.

Recall that $v$ gives charges only to neighboring $5^{-}$-vertices of $H$ by R4, and two such $5^{-}$-vertices cannot be adjacent. Among the neighbors of $v$, let $n_1$ be the number to which $v$ gives charge $1$ by R4.1, and let $n_2$ be the number to which $v$ gives charge $\rho$ by R4, where $\rho \in\{\frac{1}{5},\frac{1}{2}\}$. By the absence of $C_1$ and $C_2$, we get $n_1\leq 2$.

\noindent
\textbf{Case V4.1:} $d_G(v)=9$.

By the absence of $C_1$ and $C_4$, we have $n_1\leq1$. If $n_1=1$, then $n_2\leq 3$ and it follows that $\ch(v)\geq \ch_0(v)- 1- 3\times \frac 12>0$. Otherwise $\ch(v)\geq \ch_0(v)- 6\times \frac 12=0$.

\noindent
\textbf{Case V4.2:} $d_G(v)=10$.

If $n_1=2$, then $n_2\leq 3$ and it follows that $\ch(v)\geq \ch_0(v)- 2 \times 1- 3\times\frac 12>0$. If $n_1=1$, then $n_2\leq 5$ and so $\ch(v)\geq \ch_0(v)-1-5\times \frac 12>0$. If $n_1=0$, then $n_2\leq 6$ and we have $\ch(v)\geq \ch_0(v)- 6\times \frac 12>0$.

\noindent
\textbf{Case V4.3:} $d_G(v)\geq 11$.

If $n_1=2$, then $n_2\leq \frac{d_G(v)}{2}-2$ and $\ch(v)\geq d_G(v)- 6-2 \times 1-\frac 12(\frac{d_G(v)}{2}- 2)>0$. If $n_1=1$, then $n_2\leq \frac{d_G(v)}{2}-1$ and $\ch(v)\geq d_G(v)- 6-1-\frac 12(\frac{d_G(v)}{2}- 1)>0$. If $n_1=0$, then $n_2\leq \frac{2d_G(v)}{3}$ and $\ch(v)\geq d_G(v)- 6- \frac 12(\frac{2d_G(v)}{3})>0$.

\begin{figure}[htbp]
 \begin{center}
   \includegraphics[scale=0.75]{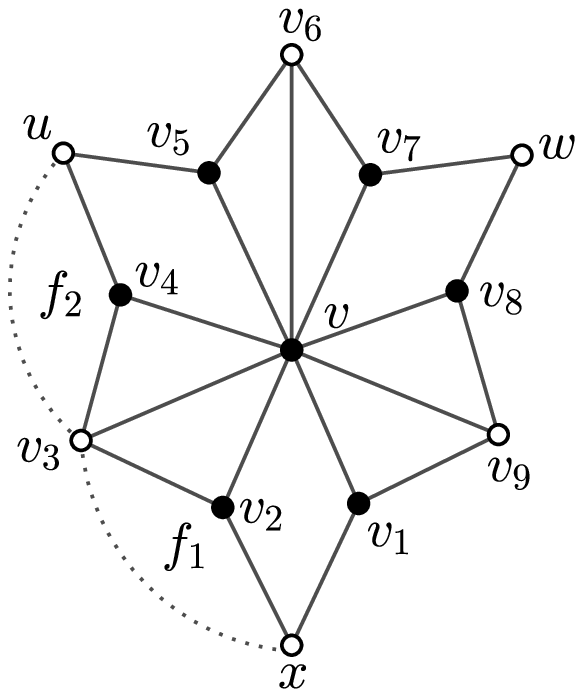}\\
   \caption{$d_G(v)=9$.}
   \label{Fig.deg9}
 \end{center}
\end{figure}

\textbf{Till now, we have checked that} \bm{$\ch(x)\geq 0$} \textbf{for every element} \bm{$x\in V(G)\cup F(G)$}.
Let $W=\{w|\, w\in V(H)\; \text{and}\; \ch(w)>0\}$. Now we prove $W\neq\emptyset$. Suppose that $W=\emptyset$, that is, $ch(v)=0$ for any $v\in V(H)$. Then $d_G(v)\leq 9$ for each $v\in V(H)$ by Case V4.2 and V4.3. If $H$ has a 9-vertex $v$, then it must be the center of the configuration shown in Figure \ref{Fig.deg9}, where every neighbor of $v$ is in $H$, for otherwise $ch(v)>0$. Let $f_1$ and $f_2$ be the faces with paths $xv_2v_3$ and $v_3v_4u$ in their boundaries, respectively. By the absence of $C_1$, $d_G(f_1)\geq 4$ or $d_G(f_2)\geq 4$, w.l.o.g. $d_G(f_1)\geq 4$. Then $v_2$ is a $3$-vertex of $H$ and is incident with only one $3$-face. Suppose that $x\in N$, so that $c(x\rightarrow v_2)=1$ by R3. It is proved in Case V2 that $ch(v_2)\geq0$, but this argument uses only $\frac{1}{2}$ of the $1$ unit of charge that $x$ gives $v_2$ by R3. Note that (ii) of Case V2.2 does not apply to $v_2$. So $ch(v_2)\geq\frac{1}{2}$ and $v_2\in W$, a contradiction. Thus $x\in V(H)$ and $d_G(x)=9$ by (b). Since $x$ is incident with two adjacent $4^{+}$-faces, it is not the center of a configuration as shown in Figure \ref{Fig.deg9}. So $\ch(x)>0$ and $x\in W$, a contradiction. Hence $d_G(v)\leq 8$ for each $v\in V(H)$.

Suppose that $H$ has a 8-vertex $v$. Since $W=\emptyset$, $ch(v)=0$. It follows from Case V3 that $v$ sends $\frac{1}{2}$ to each of four $4$-vertices by R4.2 and receives nothing. Thus all faces incident with $v$ are $3$-faces and every neighbor of $v$ is a $4$-vertex or an $8$-vertex in $V(H)$. At the same time, every adjacent $4$-vertex $u$ of $v$ is adjacent to at least three 8-vertices and receives $\frac 12$ from each of them, so $u$ is adjacent to four $8$-vertices in $V(H)$ and is incident with four 3-faces. These imply that every vertex of $H$ is a $4$-vertex or an $8$-vertex, and is not adjacent to any vertex of $N$. So $|N|=0$, a contradiction. Hence $d_G(v)\leq 7$ for each $v\in V(H)$.

By Case V3, $\ch(v)>0$ for any 7-vertex $v$ in $H$. So $d_G(v)\leq 6$ for each $v\in V(H)$. By (b), $d_G(v)=6$ for each $v\in V(H)$. Let $v$ be a vertex of $H$ that is incident with $f_0$. Then $ch(v)\geq d_G(v)-6+1+t(f_0,v)\geq 1$ by R1.2, a contradiction with $W=\emptyset$.

Hence we obtain $|W|>0$ and it follows that $\sum_{x\in{{V(G)}\cup{F(G)}}}\ch(x)>0$, a contradiction to $\mathrm{(V)}$. This completes the proof of Lemma \ref{planar}.
\end{proof}

Now we prove another property of $K_5$-minor-free graphs which is similar to that of planar graphs (see \cite{KSS,WHJ1,WWF}).

\begin{lemma}\label{Lem-k5-deg10}
Let $G$ be a connected $K_5$-minor-free graph. Suppose that
\begin{description}
    \item[$$(a)$$] $\delta(G)\geq 2$, and
    \item[$$(b)$$] $d_G(u)+d_G(v)\geq 12$ for every edge $uv\in E(G)$.
\end{description}
Then $G$ has a subgraph isomorphic to one of the configurations in Figure $\ref{fig.ve-k5-deg2}$ and $\ref{fig.ve-k5-deg3}$.
\begin{figure}[H]
 \begin{center}
   \includegraphics[scale=0.8]{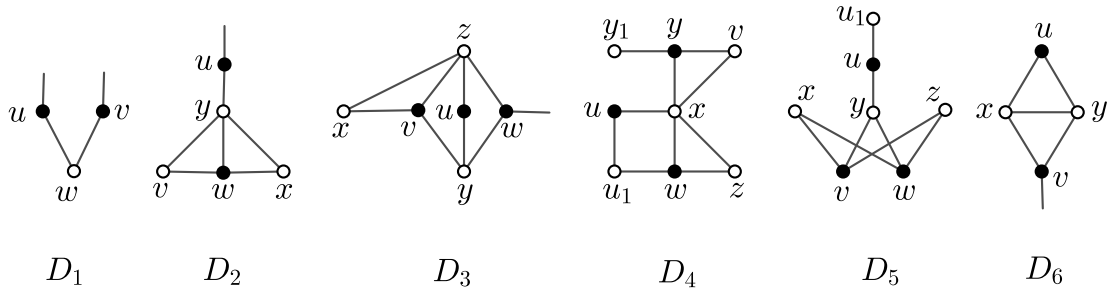}\\
   \caption{Configurations for Lemma \ref{Lem-k5-deg10}, each of which contains a $2$-vertex.}
   \label{fig.ve-k5-deg2}
 \end{center}
\end{figure}
\begin{figure}[H]
 \begin{center}
   \includegraphics[scale=0.8]{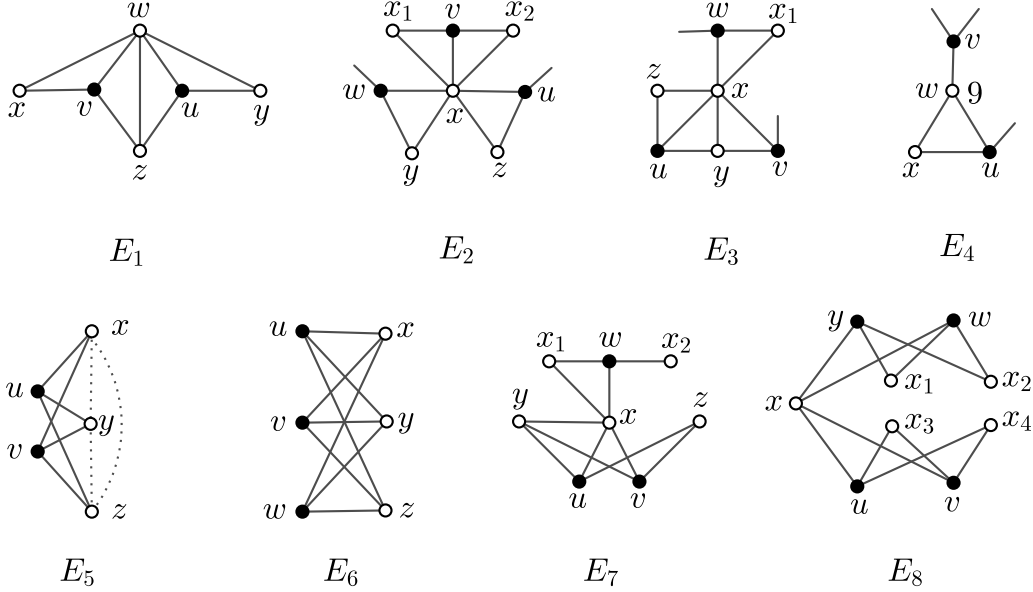}\\
   \caption{Configurations for Lemma \ref{Lem-k5-deg10}, each of which contains at least two $3$-vertices and contains no $2$-vertices. In $E_4$, $w$ is a $9$-vertex. In $E_5$, $\{xy, yz,zx\}\cap E(G)=\emptyset$. In $E_8$, $|\{x_1,x_2\}\cap \{x_3,x_4\}|\leq 1$.}
   \label{fig.ve-k5-deg3}
 \end{center}
\end{figure}
\end{lemma}
\begin{proof} Let $H$ be a counterexample to the lemma  with $|V|+|E|$ as small as possible. By the absence of $D_1$, we have that
\begin{description}
  \item[$\mathrm{(c)}$] \textit{every vertex is adjacent to at most one $2$-vertex}.
\end{description}

Let $V_2$ be the set of 2-vertices of $H$, $V_2^+=\{v\in V_2$:\, Two neighbors of $v$ are adjacent$\}$ and  $V_2^-=\{v\in V_2:$ Two neighbors of $v$ are not adjacent$\}$. Let $C=\{xy:$ there is a 2-vertex $v\in V_2^-$ such that $x,y\in N(v)\}$.
Now we construct a new graph $G$ from $H$ by letting $G=H-V_2+C$. Thus $G$ is also a $K_5$-minor free  graph with $\delta(G)\geq 3$. By (c), we have that for every $v\in V(G)$,
$$d_{G}(v)=\left\{
  \begin{array}{ll}
    d_H(v)-1, & \hbox{if $v\in \bigcup_{u\in V_2^+} N_H(u)$,} \\
    d_H(v), & \hbox{otherwise.}
  \end{array}
\right.$$
By the hypothesis (b), every neighbor of a $2$-vertex of $H$ is a $10^{+}$-vertex. So
\begin{description}
  \item[$\mathrm{(d)}$] \textit{$d_G(u)+d_G(v)\geq 12$ for each edge $uv$ of $G$.}
\end{description}
Since $\delta(G)\geq3$, $G$ does not contain any of $D_1$-$D_6$. If $G$ contains any of $E_1$-$E_8$, then that configuration must contain exactly one edge $e_0$ such that $e_0\in C$ or two ends of $e_0$ are adjacent to a 2-vertex in $V_2^+$ by (c). It is obvious that there is no such possible edge $e_0$ in $E_5$, $E_6$ or $E_8$. If $G$ contains $E_1$, then $e_0\in C$  and $e_0$ must be $wx$, $wy$ or $wz$, and $H$ contains $D_2$, $D_2$ or $D_3$ respectively, which is a contradiction. If $G$ contains $E_2$, then $e_0(\in C)$ must be $xx_1$, $xx_2$, $xy$ or $xz$, and $H$ contains $D_4$, $D_4$, $D_2$ or $D_2$ respectively. If $G$ contains $E_3$, then $e_0(\in C)$ must be $xx_1$, $xy$ or $xz$, and $H$ contains $D_2$, $D_3$ or $D_4$ respectively. If $G$ contains $E_4$, then $w,x$ must be adjacent to a 2-vertex in $V_2^+$, and $H$ contains $D_6$. If $G$ contains $E_7$, then $e_0$ must be $xx_1$ or $xy$, and $H$ contains $D_5$. In each case we have a contradiction. Hence we get

\begin{description}
  \item[$\mathrm{(e)}$] \textit{$G$ contains none of $E_1$-$E_8$.}
\end{description}

By Corollary \ref{wagner}, $G$ has a $3$-simplified tree-decomposition $(T, \mathcal{F})$ such that each part is a planar graph or the Wagner graph, and assume that $T$ contains as many vertices as possible and subject to that the number of leaves is as large as possible. By Claim \ref{cl0} and the paragraph following it, we assume $|V(T)|\geq3$ and

\begin{description}
  \item[$\mathrm{(f)}$] \textit{for each leaf $x$ of $T$,  $|S_x|=3$ and  $G^*_x=K_{1,3}$ is a star of order $4$.}
\end{description}

We choose a vertex $c$ of $T$, $N_T(c)=\{c_0, c_1, ..., c_t\}$ $(t\geq 1)$ such that $c_0$ is the parent of $c$ and $c_1, \ldots, c_t$ are the leaves of $T$. Let $K=G^*_c\cup G_{c_1}\cup\cdots\cup G_{c_t}$ and $\overline{K}=K-S_c$. By the similar arguments as Lemma \ref{Lem-k5-deg7}, we can get that $\overline{K}$ has at least one edge.

By (f), denote $V(G^*_{c_i})=\{u_i, u_{i1}, u_{i2}, u_{i3}\}$, where $u_i$ is the center of the star $G^*_{c_i}$ and $S_{c_i}=\{u_{i1}, u_{i2}, u_{i3}\}$ $(1\leq i\leq t)$. Let $G^{+}$ be described as in the proof of Lemma \ref{Lem-k5-deg7}, and note that, as in that proof, $G^{+}$ is a $K_5$-minor-free graph and $G^{+}[V_c]$ is planar. Embed $G^{+}[V_c]$ into the plane so that the vertices and edges of the complete graph $G^{+}[S_c]$ are all in the boundary of the outside face, and for each $j$ $(1\leq j\leq t)$, $u_{j1}$, $u_{j2}$, $u_{j3}$ forms some $3$-face $f_j$ of $G^{+}[V_c]$.

Suppose that $S_{c_i}\neq S_{c_j}$ for any $i, j$ $(1\leq i<j\leq t)$. Then we can embed every star $G^*_{c_j}$ into $f_j$ $(1\leq j\leq t)$ to form a new planar graph $M$. and thus $G^{+}$ is also a planar graph. Since $K\subseteq M$, $K$ is planar. Since $\delta(G)\geq3$, every vertex of $\overline{K}$ has degree at least $3$ in $K$. So $G$ contains at least one of $C_1$-$C_4$ by Lemma \ref{planar} (where $M, \overline{K}, S_c$ plays the role of $G, H, N$ in the lemma, respectively), and hence at least one of $E_1$-$E_4$ ($E_4$ is a subgraph of $C_4$), and this is the required contradiction.

So there exist $i$, $j$ $(1\leq i<j\leq t)$ such that $S_{c_i}=S_{c_j}$. If there are $i,j,k$ $(1\leq i<j<k\leq t)$ such that $S_{c_i}= S_{c_j}= S_{c_k}$, then $G$ contains $E_6$, a contradiction. Hence we can assume that

\begin{description}
  \item[$\mathrm{(g)}$] \textit{there exists an integer $s$ $(s\geq \frac{t}{2})$ such that $S_{c_i}= S_{c_{s+i}}$ for every $i$ $(1\leq i\leq t-s)$ and  $S_{c_i}\neq S_{c_j}$ for any $i$ and $j$ $(1\leq i<j\leq s)$.}
\end{description}

Let $K^{+}$ be a planar graph from $G^{+}[V_c]$ by embedding every star $G^*_{c_j}$ into $f_j$ $(1\leq j\leq s)$ and let $K'=G^*_c\cup G_{c_{1}}\cup\cdots\cup G_{c_s}$. Then $K'= K-\{u_{s+1}, u_2, \cdots, u_t\}$, $K'\subseteq K^{+}$ and $K'$ is a planar graph. Since $G$ has no $E_8$, $S_{c_i}\cap S_{c_j}=\emptyset$ whenever $s+1\leq i<j\leq t$. So for every $v\in V(K')$,
$$d_{K'}(v)=\left\{
  \begin{array}{ll}
    d_K(v)-1, & \hbox{if $x\in \bigcup^t_{i=s+1}S_{c_i}$,} \\
    d_K(v), & \hbox{otherwise.}
  \end{array}
\right.\eqno{\mathrm{(VI)}}$$

Let $K''=K'-S_c$. For a vertex $v\in V(K'')$, $d_K(v)=d_G(v)\geq3$ and $d_{K^{'}}(v)=d_K(v)$ if $d_K(v)\leq8$, so we have $d_{K^{'}}(v)\geq3$. Since $N(u_i)=N(u_{i+s})$ and $\overline{K}$ has at least one edge, $K''$ has at least one edge. If there is an edge $xy\in E(K'')$ such that $d_{K'}(x)\leq d_{K'}(y)$ and $d_{K'}(x)+ d_{K'}(y)<12$, then there exits a $i$ $(1\leq i \leq t-s)$ such that $x\in u_i$ and $y\in S_{c_i}$, w.l.o.g. assume that $y=u_{i1}$. Then $d_{K'}(u_i)=3$, $d_{K'}(u_{i1})=8$ and it follows that $d_{K}(u_{i1})= 9$, $u_{i1}u_{i2}, u_{i1}u_{i3}\not\in E(K)$ by the absence of $E_4$, and $u_{i2}u_{i3}\in E(K)$ by the absence of $E_5$. This implies that $d_{K'}(u_{i2})\geq9$ and $d_{K'}(u_{i3}) \geq 9$.

So we may assume that every edge $xy\in E(K'')$ such that $d_{K'}(x)+d_{K'}(y)<12$ is of the form $u_iu_{i1}$ for some $i$ $(1\leq i\leq t-s)$. For every such edge, add the edge $u_{i1}u_{i2}$ to $K'$ to form a new graph $K^{*}$, which is planar since $K^{*}\subseteq K^{+}$. Let $H^{*}=K^{*}-S_c$. Then $d_{K^{*}}(x)+d_{K^{*}}(y)\geq12$ for every edge $xy\in E(H^{*})$, $d_{K^{*}}(v)\geq d_{K^{'}}(v)\geq3$ for every vertex $v\in V(H^{*})$ and $H^{*}$ has at least one edge. Hence $K^{*}$ contains at least one of $C_1$-$C_4$ by Lemma \ref{planar}.

Suppose that $K^*$ contains $C_2$. We label the vertices of $C_2$ as $E_2$. Since $G$ does not contain $E_2$, there is an edge $u'v'\in \{xx_1, xx_2, xy, xz\}$ such that $u'v' \in E(K^*)\backslash E(K)$, that is, $u'v'$ is a new added edge as above from $K'$ to $K^*$. Suppose that $u'v'=xz$ (we can settle the case $u'v'=xy$ similarly), w.l.o.g. assume that $u'=x=u_{s+1,1}$ and $v'=z=u_{s+1,2}$. Then $d_{K'}(x)=8$. So $d_{G}(x)=9$ or $d_{G}(z)=9$. If $G[N_G(u)]$ is an empty graph, then $G$ contains $E_5$. Otherwise, $G$ contains $E_4$. If $u'v'\in \{xx_1,xx_2\}$, then w.l.o.g. $v=u_{s+1}$, $N(u_1)=N(v)$, and it follows that $G$ contains $E_7$. By the similar arguments, if $K^*$ contains $C_3$, then $G$ contains $E_4$, $E_5$ or $E_7$. If $K^*$ contains $C_1$ or $C_4$, then $G$ contains $E_7$. These contradictions complete the proof of Lemma \ref{Lem-k5-deg10}.
\end{proof}

\section{Proof of Theorem \ref{th1}} \label{section4}

The coloring methods of $(1)$-$(3)$ in Lemma \ref{Lem-k5-deg7} can be found in \cite{Bo}, and all configurations in Figure \ref{Fig.planar7} are reducible (see \cite{SZ}). So Theorem \ref{th1}(1) holds.

Next we give the proof of Theorem \ref{th1}(2) by using Lemma \ref{Lem-k5-deg10}. To complete the proof, it suffices to prove that if an integer $k\geq  11$ and $G$ is $K_5$-minor-free graph $G$ with $\Delta(G)\leq k-1$, then $\chi''(G)\leq k$.

A \emph{nice} $k$-coloring, or simply nice coloring, of $G$ is a $k$-coloring of all its edges and all its vertices of degree at least $6$ such that no two adjacent or incident elements receive the same color. The uncolored vertices are then easily colored, as each is adjacent to at most five vertices and incident with at most five edges, and so there is a color available for it. Thus $G$ has a total $k$-coloring if and only if it has a nice coloring. We will work with nice colorings throughout the remaining parts.

Given a partial nice coloring of a graph, we say that a color is \emph{present} at a vertex $v$ if it used on $v$ or on an edge incident with $v$, otherwise it is \emph{missing} at $v$. Let $G$ be a minimal $K_5$-minor-free graph such that $\Delta(G)<k$ and $G$ has no nice coloring, i.e., no total $k$-coloring. We will get a contradiction.

If $G$ contains a $1$-vertex $v$, consider a nice coloring of $G-v$. The edge at $v$ is incident with one colored vertex and adjacent to at most $\Delta(G)-1$ colored edges, and so there is a color available for it. This gives a nice coloring of $G$, which contradicts the hypothesis of $G$.

In a similar way, if $G$ contains an edge $uv$ such that $d(u)+d(v)\leq11$, consider a nice coloring of $G-uv$. Then $uv$ is incident with at most one colored vertex and $d(u)+d(v)-2\leq9$ colored edges, and so it can be colored so as to give a nice coloring of $G$.

These contradictions show that $\delta(G)\geq2$ and that $G$ does not have an edge $uv$ as just described. It now follows from Lemma \ref{Lem-k5-deg10} that $G$ must contain one of the fourteen configurations $D_1$-$D_6$ and $E_1$-$E_8$. In the following, we will show that all of these configurations are reducible, that is, they cannot occur in our minimal counterexample $G$.

Many of these configurations have already been proved reducible elsewhere (for specific value of $\Delta$, but the proofs work for general $\Delta$), as follows.
\begin{table}[H]
\small
\centering
\begin{tabular}{lll}
$D_1$: \cite{KSS}, Lemma 2(iii)&$D_6$: \cite{KSS}, Lemma 3(iv)&$E_3$: \cite{Hou}, Fig. 1(b)\\
$D_2$: \cite{KSS}, Fig. 2(b) &$E_1$: \cite{KSS}, Fig. 5(a)&$E_4$: \cite{KSS}, Fig. 6(a) \\
$D_3$: \cite{KSS}, Fig. 2(a) &$E_2$: \cite{KSS}, Fig. 3(h)& \\

\end{tabular}
\label{tab1}
\end{table}

It remains to prove the reducibility of configurations $D_4$, $D_5$, $E_5$, $E_6$, $E_7$ and $E_8$. We start with $E_5$ and $E_6$, which are the easiest.

\noindent
\textbf{Configuration $E_5$.} Let $\psi$ be a nice coloring of $G-\{u,v\}+\{xy, yz, xz\}$ by the minimality of $G$. We color $ux, vy$ with $\psi(xy)$, $uy, vz$ with $\psi(yz)$, and $uz, vx$ with $\psi(zx)$, to give a nice coloring of $G$, a contradiction.

\noindent
\textbf{Configuration $E_6$.} The edges at $u$, $v$ and $w$ form a $K_{3,3}$. Given a nice coloring of $G-\{u,v,w\}$, each edge of this $K_{3,3}$ has at least three available colors, and so these edges can be colored by the result of Galvin \cite{FG} that for a bipartite graph, list edge chromatic number equals to its maximum degree. This gives a nice coloring of $G$, which is a contradiction.

In the proofs for configurations $D_4$ and $D_5$, we will use the following simple claim.
\begin{claim}\label{cl1}
Suppose $u$ is a $2$-vertex of $G$ with neighbors $x$ and $y$, and $\psi$ is a nice coloring of $G-uy$ such that $\psi(ux)=a$. Then \\
$(i)$ $a$ is the unique color missing at $y$, and \\
$(ii)$ $a$ is present at every neighbor $z$ of $y$ $($on an edge at $z$ if $d(z)\leq5$, since then $z$ is uncolored$)$.
\end{claim}
\noindent
\textbf{Proof of Claim \ref{cl1}.} At least one color is missing at $y$, and if $y$ misses $b\neq a$, then we can color $uy$ with $b$ to give a nice coloring of $G$. This contradiction proves $(i)$. If $a$ is missing at $z$ then we can recolor $yz$ with $a$, which contradicts $(i)$. This proves $(ii)$. \ \ \ \ $\blacksquare$

\noindent
\textbf{Configuration $D_4$.} Let $\psi$ be a nice coloring of $G-ux$ and let $a=\psi(uu_1)$ and $b=\psi(u_1w)$. By Claim \ref{cl1}, $a$ is the unique color missing at $x$, and $a$ is present at $w$ and $y$, which forces $\psi(wz)=a$. If $\psi(xz)\neq b$, then we could swap the colors of $xz$ and $wz$ so that $a$ is present at $x$, which contradicts Claim \ref{cl1}. So $\psi(xz)=b$. If $b$ is not present at $y$, then we could swap colors $a$ and $b$ along the path $uu_1wzx$, recolor $xy$ with $b$, color $ux$ with $\phi(xy)$ to obtain a nice coloring of $G$, a contradiction. So $b$ is present at $y$, which means that the colors of $vy$ and $yy_1$ are $a$ and $b$ in some order. Now swap the colors of $vy$ and $vx$ and color $ux$ with $\psi(vx)$ firstly. Then if $\psi(vy)=b$,  we swap colors $a$ and $b$ along the path $uu_1wzx$. Thus we also obtain a nice coloring of $G$, a contradiction.

\noindent
\textbf{Configuration $D_5$.} Let $\psi$ be a nice coloring of $G-uy$ and let $a=\psi(uu_1)$. By Claim \ref{cl1}, $a$ is the unique color missing at $y$, and $a$ is present at $v$ and $w$, w.l.o.g. $\psi(vx)=\psi(wz)=a$.

Since configuration $E_5$ is reducible, $\{xy, xz, yz\}\cap E(G)\neq \emptyset$. If $xz\in E(G)$, assume by symmetry that $\psi(xz)\neq \psi(vy)$, then we swap colors along the path $vxzw$ and recolor $wy$ with $a$ to obtain that $a$ is present at $y$, which contradicts Claim \ref{cl1}.
Suppose that $xy\in E(G)$(we can settle the case $xy\in E(G)$ similarly). If $\psi(vz)=\psi(xy)= b$ then we swap the colors $a$ and $b$ along the path $yxvzw$, otherwise we just swap the colors of $yx$ and $xv$. In either case $a$ is now present at $y$, a contradiction.

The following claim will be used  to prove Configuration $E_7$ and $E_8$.
\begin{claim}\label{cl2}
Suppose that $G$ contains two $3$-vertices $u$, $v$ with the same three neighbors $x$, $y$, $z$. Let $\psi$ be a nice coloring of $G-\{u,v\}$ and $C_t$ be a set of two colors missing at $t$ for any $t\in \{x,y,z\}$. Then
\begin{description}
  \item[$(i)$] $C_x\cap C_y\cap C_z\neq \emptyset$;
  \item[$(ii)$] If $xy\in E(G)$, $C_x=\{1,2\}$, $\psi(xy)=3$ and $1\in C_x\cap C_y\cap C_z$, then either $C_z=\{1,3\}$ or $C_x=C_y=C_z=\{1,2\}$;
  \item[$(iii)$] With the colors as in $(ii)$, let $w$ be a neighbor of $x$, $w\notin\{u,v,y,z\}$. Then color $1$ is present at $w$.
\end{description}
\end{claim}

\noindent
\textbf{Proof of Claim \ref{cl2}.} $(i)$ Suppose first that $C_x\cap C_y=\emptyset$. Color $uz$, $vz$ with the colors in $C_z$. There are two ways of coloring $ux$, $vx$ from the colors in $C_x$, and at least one of them does not clash with the colors of $uz$, $vz$. The same is true of coloring $uy$, $vy$ with the colors in $C_y$. Thus there is a nice coloring of $G$, which is a contradiction. So we may assume that no two of $C_x$, $C_y$, $C_z$ are disjoint. If $(i)$ is false, then these sets must be of the form $\{1,2\}$, $\{1,3\}$, $\{2,3\}$, but then a nice coloring of $G$ is easily completed. (Color one of the uncolored edges at random, and the colors of the rest are then uniquely determined.) Thus $(i)$ holds.

\vspace{3mm}
\noindent
$(ii)$  Suppose that $C_z=\{1,\alpha\}\neq \{1,3\}$.  Recolor $xy$ with $1$, this changes $C_x$ from $\{1,2\}$ to $\{2,3\}$ without changing $C_z$. For $(i)$ to still hold it must be that $\alpha=2$ and it follows that $C_y=\{1,2\}$. So  $C_x=C_y=C_z=\{1,2\}$. This proves $(ii)$.

\vspace{3mm}
\noindent
$(iii)$ Suppose color $1$ is missing at $w$, and let $\psi(xw)=4$. Recolor $xw$ with $1$, this changes $C_x$ from $\{1,2\}$ to $\{2,4\}$ without changing $C_z$. If $C_z=\{1,3\}$, then we have an immediate contradiction with $(i)$, otherwise $C_z=\{1,2\}\neq C_x$, which contradicts $(ii)$. \hfill $\blacksquare$

\vspace{3mm}
Now we begin to prove that configurations $E_7$ and $E_8$ are reducible.

\vspace{2mm}
\noindent
\textbf{Configuration $E_7$.} Let $\psi$ be a nice coloring of $G-\{u,v\}$ with the colors and definitions from Claim \ref{cl2}$(ii)$, and let $\psi(xx_1)=4$. By Claim \ref{cl2}$(iii)$, color $1$ is present at $w$. Suppose first that $C_x=C_y=C_z=\{1,2\}$. Then color $2$ is present at $w$ by symmetry, and so $\{\psi(wx_1), \psi(wx_2)\}=\{1,2\}$. Swapping the colors of $wx_1$ and $xx_1$ makes $C_x=\{1,4\}$ or $\{2,4\}$ while leaving $C_z=\{1,2\}\neq C_x$ and $\psi(xy)=3$. This contradicts Claim \ref{cl2}$(ii)$.

So we may assume that $C_z=\{1,3\}$. Recolor $xy$ with $1$, so that $3\in C_x\cap C_y\cap C_z$. By Claim \ref{cl2}$(iii)$, we see that $\{\psi(wx_1), \psi(wx_2)\}=\{1,3\}$. Swap the colors of $wx_1$ and $xx_1$, so that $xx_1$ has color $1$ or $3$, and give $xy$ the opposite color. This changes $C_x$ from $\{1,2\}$ to $\{2,4\}$ while leaving $C_z=\{1,3\}$, and this contradicts Claim \ref{cl2}$(i)$.

\vspace{2mm}
\noindent
\textbf{Configuration $E_8$.} Since the configuration $E_5$ is reducible, neither $\{x,x_1,x_2\}$ nor $\{x,x_3,x_4\}$ is an independent set of $G$. There are two cases.

\vspace{2mm}
\noindent
\textbf{Case 1.} $\{xx_1, xx_2, xx_3, xx_4\}\cap E(G)=\emptyset$, and so $x_1x_2$ and $x_3x_4\in E(G)$.

\vspace{2mm}
Since configurations $E_5$ and $E_7$ are reducible, $\{x_1,x_2\}\cap\{x_3,x_4\}=\emptyset$. Let $\psi$ be a nice coloring of $G'=G-\{u,v,w,y\}+\{xx_1, xx_2, xx_3, xx_4\}$. W.l.o.g., assume that $\psi(xx_i)=i$ ($i\in\{1,2,3,4\}$). For $i=1,2,3,4$, let $C_{x_i}$ be a set of two colors consisting of $i$ and another color missing at $x_i$ in $G$. We may assume that $\psi(x_1x_2)\neq3$ (for otherwise, we can swap colors $3$ and $4$ throughout $G$, and swap the labels $x_3$ and $x_4$). In a similar way, we may assume that $\psi(x_3x_4)\neq2$.

Let $C_x=\{1,3\}$. Note that $\psi(x_1x_2)\notin C_x$, and $C_x\neq C_{x_2}$ since $2\in C_{x_2}$. It follows from the proof of Claim \ref{cl2}$(i)$ and $(ii)$ (with $x_2$, $x_1$, $x$ in place of $x$, $y$, $z$ respectively) that, after recoloring $x_1x_2$ if necessary, we can color the six edges at $w$ and $y$ with the edges $xw$ and $xy$ having colors $1$ and $3$. By the same argument, we can color the six edges at $u$ and $v$ with edges $xu$ and $xv$ having colors $4$ and $2$. This gives a nice coloring of $G$, a contradiction.

\vspace{2mm}
\noindent
\textbf{Case 2.} $\{xx_1, xx_2, xx_3, xx_4\}\cap E(G)\neq\emptyset$, w.l.o.g. $xx_3\in E(G)$.

\vspace{2mm}
\noindent
\textbf{Subcase 2.1.} $\{x_1,x_2\}\cap\{x_3,x_4\}=\emptyset$ or $x_3\notin\{x_1,x_2\}$.

\vspace{2mm}
Since the configuration $E_7$ is reducible, $x$ cannot be adjacent to $x_1$ or $x_2$, and so $x_1x_2\in E(G)$. Given a nice coloring $\psi$ of $G-\{u,v\}$, let $C_x$, $C_{x_3}$ and $C_{x_4}$ be sets of two colors missing at $x$, $x_3$ and $x_4$ respectively, and by Claim \ref{cl2} assume that $1\in C\cap C_{x_3}\cap C_{x_4}$, $C_x=\{1,2\}$, $\psi(xx_3)=3$, $\psi(xy)=4$, $\psi(xw)=5$ and $C_{x_4}=\{1,2\}$ or $\{1,3\}$. By Claim \ref{cl2}$(iii)$, color $1$ is present at $w$ and $y$, say $\psi(wx_1)=\psi(yx_2)=1$. Swap colors along the path $wx_1x_2y$, if $\psi(x_1x_2)=4$ then recolor $xy$ with $1$, otherwise recolor $xw$ with $1$. This changes $C_x$ from $\{1,2\}$ to $\{2,4\}$ or $\{2,5\}$. If $C_{x_4}=\{1,3\}$, this contradicts Claim \ref{cl2}$(i)$. Otherwise it contradicts Claim \ref{cl2}$(ii)$.

\vspace{2mm}
\noindent
\textbf{Subcase 2.2.} $\{x_1,x_2\}\cap\{x_3,x_4\}\neq\emptyset$ and $x_3\in\{x_1,x_2\}$, w.l.o.g. assume that $x_1=x_3$.

\vspace{2mm}
Then we have $x_2\neq x_4$ since $|\{x_1,x_2\}\cap \{x_3,x_4\}|\leq1$ (see the caption of Figure \ref{fig.ve-k5-deg3}), and $x_1x_2, x_3x_4\notin E(G)$ since the configuration $E_7$ is reducible.  Given a nice coloring $\psi$ of $G-\{u,v\}+\{xx_4, x_3x_4\}$. Let $\psi(x_3x_4)=1$, $\psi(xx_4)=2$, $\psi(xx_3)=3$, $C_{x_4}=\{1,2\}$, $C_{x_3}=\{1,\alpha\}$ and $C_{x}=\{2,\beta\}$ where $\alpha$ and $\beta$ are colors missing at $x_3$ and $x$ respectively. By Claim \ref{cl2}$(i)$ and $(ii)$, $\beta=1$ and $\alpha=2$. So $C_{x_4}=C_{x_3}=C_{x}=\{1,2\}$. It implies that $\{\psi(xw), \psi(xy), \psi(x_3w), \psi(x_3y)\}\cap\{1,2\}=\emptyset$. It follows from Claim \ref{cl2}$(iii)$
that colors $1$ and $2$ are present at $w$ and $y$, and this is the required contradiction.

Hence we complete the proof of Theorem \ref{th1}.

\section{Conclusion}
Here, we proved that if $G$ is a $K_5$-minor-free graph and $\Delta(G)\geq 10$, then $\chi''(G)=\Delta(G)+1$. But Kowalik et al. \cite{KSS} proved that the same result about planar graphs with $\Delta(G)\geq 9$, there is still a gap in terms of maximum degree constraints. We think that it is also true for $K_5$-minor-free graphs and pose a similar conjecture as in \cite{Shen}.

\begin{conjecture}
If $G$ is a $K_5$-minor-free graph and $\Delta(G)\geq 4$, then $\chi''(G)=\Delta(G)+1$.
\end{conjecture}

Finally, it is worth mentioning that our results also hold for $K_{3,3}$-minor-free graphs. Before proceeding, we present another remarkable result proved by Wagner in \cite{Wagner}.
\begin{lemma}{\rm \cite{Wagner}}
Let $G$ be an edge-maximal graph without a $K_{3,3}$ minor. If $|G|\geq 4$ then $G$ has a $2$-simplicial tree-decomposition $(T, \mathcal{F})$ such that each part is a planar triangulation or a copy of the complete graph $K_5$.
\end{lemma}

According to the lemma, it is similar to prove that any $K_{3,3}$-minor-free graph satisfies Lemma \ref{Lem-k5-deg7} and Lemma \ref{Lem-k5-deg10}. Thus we can obtain the following theorem.

\begin{theorem}
Let $G$ be a $K_{3,3}$-minor-free graph. Then $\chi''(G)\le \max\{9, \Delta(G)+2\}$. Moreover, $\chi''(G)=\Delta(G)+1$ if $\Delta(G)\geq 10$.
\end{theorem}

Finally, we note that Feng, Gao and Wu \cite{FGW} obtained a result about the edge coloring of $K_5$-minor-free graphs. Very recently, Feng and Wu and Yang \cite{FWY} consider a list edge coloring of $K_5$-minor-free graphs. It is worth to investigate the list total coloring of $K_5$ (or $K_{3,3}$)-minor-free graphs.

\noindent
\textbf{Acknowledgments}

We are very graceful to the referees for their careful reviewing and valuable suggestions. It is worth mentioning that they provide us with new statements in the proof of Lemma \ref{planar-deg7}, Lemma \ref{Lem-k5-deg7}, Lemma \ref{planar} and Theorem \ref{th1}(2), which greatly improve this paper.

\end{document}